\documentclass[12pt,a4paper,reqno,twoside]{amsart}

\usepackage[english]{babel}
\usepackage{stmaryrd}
\usepackage{dsfont}
\usepackage{bbm}
\usepackage[symbol*,ragged]{footmisc}
\usepackage[colorlinks,linkcolor=red,anchorcolor=blue,citecolor=blue,urlcolor=blue]{hyperref}
\usepackage{color,xcolor}

\usepackage{geometry}
\usepackage{amssymb}
\usepackage{amsmath}
\usepackage{mathrsfs}
\usepackage{amsfonts}
\usepackage{epsfig}

\usepackage{amsthm}
\usepackage{amsxtra}
\usepackage{bbding}
\usepackage{epsfig}
\usepackage{graphicx}
\usepackage{latexsym}
\usepackage{mathbbol}
\usepackage{bbold}

\usepackage{pifont}
\usepackage{wasysym}
\usepackage{skull}
\usepackage{float}

\DeclareSymbolFontAlphabet{\mathbb}{AMSb}
\DeclareSymbolFontAlphabet{\mathbbol}{bbold}

\usepackage{amscd}
\usepackage[all]{xy}
\allowdisplaybreaks[4]
\usepackage{setspace}

\theoremstyle{plain}

\newtheorem*{corollary*}{\normalfont\scshape Corollary}

\theoremstyle{remark}
\newtheorem*{remark*}{\normalfont\scshape Remark}

\numberwithin{equation}{section}
\addtocounter{footnote}{1}

\renewcommand{\footnoterule}{
  \kern -3pt
  \hrule width 2.5in height 0.4pt
  \kern 3pt
}

\makeatletter
\@ifundefined{MakeUppercase}{}{}
\makeatother

\begin{document}
	
\title{On  some sums involving small arithmetic functions}

\author {Wenguang Zhai}

\address{[Wenguang Zhai]  Department of Mathematics, China University of Mining and Technology,          Beijing 100083, People's Republic of China}

\email{zhaiwg@hotmail.com}

\date{}


\thanks{This work is  supported by the National Natural Science Foundation of China( Grant No. 11971476). The author deeply thanks the referees for careful reading of the manuscript and for many valuable suggestions.  \\
$\clubsuit$ This paper was published in {\it Acta Mathematica Sinica(English Series)},  2024, Vol. {\bf 40}, No. 10, pp. 2497-2518, which was
   accepted for publication in November 24, 2022.}

\begin{abstract}
  Let $f$ be any arithmetic function and define $S_f(x):=\sum_{n\leq x}f([x/n])$.
  If the function $f$ is small, namely, $f(n)\ll n^\varepsilon,$ then the error term $E_f(x)$ in the asymptotic formula of  $S_f(x)$ has the form
  $O(x^{1/2+\varepsilon}).$  In this paper, we shall study the mean square of $E_f(x)$ and  establish some new results of $E_f(x)$ for
  some special functions. 
\end{abstract}

\maketitle

\addtocounter{footnote}{1}

\section{Introduction}

\subsection{On a special sum }

Let $f: {\Bbb N}\rightarrow {\Bbb C}$ be any arithmetic function. For any $x\geq 1,$ define
\begin{equation}
S_f(x):=\sum_{n\leq x}f ( [ x/n ] ).
\end{equation}

This sum was first studied in  \cite{BHS} and then was extensively studied by many authors.
 See,  for example,  \cite{ Bo, Hey,  LWY, LWY2, MS0,MS, MW, MRWJ, Sj, Wu, Wu2, Zh,Zh2}.

Wu \cite{Wu} and Zhai \cite{Zh} proved independently that
if $f(n)\ll n^\alpha(\log n)^\theta$ for some $0\leq \alpha<1$ and $\theta\geq 0,$ then the asymptotic formula
\begin{equation}
S_f(x)=C_f x+O(x^{\frac{1+\alpha}{2}}(\log x)^\theta)
\end{equation}
holds, where the constant  $C_f$ is defined by
\begin{equation}
C_f:=\sum_{n=1}^\infty \frac{f(n)}{n(n+1)}.
\end{equation}

If $f(n)\ll n^\varepsilon$, then we say $f$ is a {\it small} arithmetic function. From (1.2) we see that if $f$ is any {\it small}  arithmetic function, then
\begin{equation}
S_f(x)=C_f x+O(x^{1/2+\varepsilon} ).
\end{equation}

The exponent $1/2$ in the error term of (1.4) is a barrier for small functions $f$'s.
 It is interesting and natural to ask the following question:
 for a given small arithmetic function $f,$ is it possible to break the barrier $1/2?$

 Several special examples have been studied. For simplicity, define $$E_f(x):=S_f(x)-C_f x$$ and
 let $\theta_f$ denote the infimum of $\alpha_f$ for which the estimate $E_f(x)\ll x^{\alpha_f+\varepsilon}$ holds.
 Ma and Wu \cite{MW} first proved that $\theta_\Lambda\leq 35/71$, where $\Lambda$ is the von Mangoldt function.
The exponent $35/71$ was improved to $97/203$ and $9/19 $ in Bordell\'{e}s \cite{Bo},  Liu, Wu and Yang\cite{LWY}, independently.
Ma and Sun  \cite{MS} proved that if  $f(n)=\tau(n),$ the Dirichlet divisor function, then
 $\theta_\tau\leq 11/23.$ The exponent $11/23$ was improved to $19/40$ and $ 5/11$ in Bordell\'{e}s \cite{Bo} and Stucky \cite{Sj}, independently.

Bordell\'{e}s \cite{Bo} also studied some other arithmetic functions. He proved that
\begin{eqnarray}
&&\theta_{\tau_3}\leq \frac{283}{574},\ \ \theta_{\tau_k}\leq \frac 12-\frac{1}{2(4k^3-k-1)}(k\geq 4),\  \theta_{\omega}\leq \frac{455}{914},\ \
\theta_{2^\omega}\leq \frac{97}{202},
\end{eqnarray}
where $\tau_k(n)$ denotes the number of ways $n$ can be written as a product of $k$ factors, $\omega(n)$ denotes the number of distinct prime
divisors of $n$, respectively. All results in (1.5) were
improved by Liu, Wu and Yang \cite{LWY2}, where they proved that
\begin{equation}
  \theta_{\tau_k}\leq   \frac{5k-1}{10k-1}(k\geq 3),\ \ \theta_{\omega}\leq \frac{53}{110},\ \
\theta_{2^\omega}\leq \frac{9}{19}.
\end{equation}

Recently, Zhai \cite {Zh2} proved that if a small function $f$ satisfies a good {\it binary additive} property, then
the barrier $1/2$ in the asymptotic formula (1.4) can be broken.

\subsection{\bf Some new results of $S_f(x)$}

We first study the mean square of $E_f(x)$ for arbitrary $f$. We have the following Theorem 1.

{\bf Theorem 1.}  {\it  Let $f$ be any small arithmetic function and  $T\geq 10$ be a large parameter. Then we have
\begin{equation}
\int_T^{2T}|E_f(x)|^2dx\ll_{f,\varepsilon} T^{9/5+\varepsilon}.
\end{equation}
}

{\bf Remark 1.} From (1.4) we have $E_f(x)\ll x^{1/2+\varepsilon}$.   Theorem 1  implies  that the estimate $E_f(x)\ll x^{2/5+\varepsilon}$ holds on average.

From Theorem 1 we propose the following conjecture.

{\bf Conjecture 1.}   {\it Let $f$ be any small arithmetic function. Then we have
$$
S_f(x)=C_f x+O(x^{2/5+\varepsilon}).
$$ }

Now we study some special small arithmetic functions such that we can break the barrier $1/2$.
For any arithmetic function $f(n),$ we define
\begin{equation}
f_\ell(n):=\sum_{n=n_1\cdots n_\ell} f(n_1)\cdots f(n_\ell),\ \ f_1(n)=f(n).
\end{equation}

{\bf Theorem 2.}  {\it   Let $k\geq 2$ and $\ell\geq 1$ be two  fixed integers, and
 suppose $f$ is any function in the set  $ \{\tau_k, \Lambda_\ell, \mu_\ell, \omega_\ell\},$
where $\mu(\cdot)$ denotes the M\"{o}bius function.
Then   we have  the asymptotic formula
$$
S_f(x)=C_f x+  O\left(x^{8/17+\varepsilon}\right).
$$
}

Let ${\Bbb P}$  denote  the set of all prime numbers and let ${\bf 1}_{\Bbb P}$ denote its characteristic function.
Heyman \cite{Hey} proved that
\begin{equation}
S_{{\bf 1}_{\Bbb P}}(x):=\sum_{n\leq x}{\bf 1}_{\Bbb P}\left(\left[\frac xn\right]\right)=C_{{\bf 1}_{\Bbb P}}x+O(x^{1/2}).
\end{equation}
R. Ma and J. Wu \cite{MRWJ} proved that the exponent $1/2$ in (1.9) can be replaced by $9/19+\varepsilon.$
By the arguments of \cite{MRWJ} and Theorem 2 with the case $f=\Lambda_1$
we get the following Corollary 1.

{\bf Corollary 1.}  {\it We have the asymptotic formula
\begin{equation}
S_{{\bf 1}_{\Bbb P}}(x) =C_{{\bf 1}_{\Bbb P}}x+O(x^{8/17+\varepsilon}).
\end{equation}

}

{\bf Remark 2.} Theorem 2 is  uniform
for $ k\geq 2$ and $\ell\geq 1,$ which improves all results listed in (1.6). Although our result for
$\tau_2$ is weaker than that of Stucky \cite{Sj}, we keep it for completeness and the following Theorem 3.

{\bf Theorem 3.}  {\it   Let $k\geq 2$ and $\ell\geq 1$ be two  fixed integers, and suppose $f$ is any function in the set  $ \{\tau_k, \Lambda_\ell, \mu_\ell, \omega_\ell\}.$   Suppose $g: {\Bbb N}\rightarrow {\Bbb C}$  is any small arithmetic function and define
$F(n):=\sum_{n=d^2m}f(m)g(d).$ We have  the asymptotic formula
$$ S_F(x)=C_F x+   O\left(x^{73/155+\varepsilon}\right).$$
 }

It is well-known that
$$2^{\omega(n)}=\sum_{n=d^2m}\tau(m)\mu(d),\ \ \tau(n^2)=\sum_{n=md^2}\tau_3(m)\mu(d),\ \ \tau^2(n)=\sum_{n=md^2}\tau_4(m)\mu(d).$$
The Liouville function $\lambda(n)$ satisfies
\begin{eqnarray*}
  \lambda(n)=(-1)^{\Omega(n)}=\sum_{n=md^2}\mu(m),
\end{eqnarray*}
where $\Omega(m)$ denotes the total number of   prime divisors of $n.$
 Define the divisor functions $$t_1(n):=\sum_{n=n_1n_2n_3^2}1,\  \ t_2(n):=\sum_{n=n_1n_2n_3^2n_4^2}1.$$
 These two functions $t_1(n)$ and  $t_2(n)$ are important
when counting subgroups of finite abelian groups; see, for example, \cite{Wu3,Zh3,ZC}.

  From Theorem 3 we get the following

{\bf Corollary 2.}  {\it Let $f(n)\in \{2^{\omega(n)},\tau(n^2), \tau^2(n), \lambda(n),  t_1(n), t_2(n)\}$ and  $\ell\geq 1$ be a fixed integer. Then   we have the asymptotic formula
$$
S_{f_\ell}(x)=C_{f_\ell} x+  O\left(x^{73/155+\varepsilon}\right),
$$
where $f_\ell$ was defined in (1.8).
}

{\bf Remark 3.} From the proofs of Theorem 2 and Theorem 3 we see that if   $ f_1, f_2\in  \bigcup_{\ell\geq 1}\{\tau_{\ell+1}, \Lambda_\ell, \mu_\ell, \omega_\ell\} ,$ then Theorem 2 and Theorem 3 also hold for $f(n)=f_1\ast f_2(n)=\sum_{n=n_1n_2}f_1(n_1)f_2(n_2).$

The structure of this paper is as follows. In Section 2 we give some lemmas needed for the proofs.
In Section 3 we give the proof of Theorem 1. In Section 4 we give some estimates of exponential sums,
which are important for the proofs  of our theorems.
We
give the proofs of Theorem 2 and Theorem 3 in Section 5 and Section 6, respectively.

{\bf Notation.} Throughout this paper, $\tau_k(n)$ denotes the general divisor function, which counts
the number of ways $n$ can be written as a product of $k$ factors, $\tau_2(n)=\tau(n),$
$\Lambda(n)$ denotes the von Mangoldt function,
$\mu(n)$ denotes the
 M\"{o}bius function, $\omega(n)$ denotes the number of distinct prime
divisors of $n$, $\Omega(m)$ denotes the total number of   prime divisors of $n,$  respectively. We use
${\Bbb N}$ and ${\Bbb C}$  denote the set of positive integers and the set of complex numbers, respectively.
  For a real number $t,$ $[t]$ denotes its integer part, $\{t\}$ denotes its fractional part,
$\psi(t)=\{t\}-1/2,$  $\Vert t\Vert=\min(\{t\},1-\{t\})$ and $e(t)=exp(2\pi it).$
The symbol $n\sim N$ means that the summation condition of $n$ is $N<n\leq 2N$ and $n\asymp N$ means that
$c_1N\leq n\leq c_2N$ for two absolute positive constants $0<c_1<c_2$.
We always use $\varepsilon$ to denote a small positive constant, which maybe different
at different places. In Lemma 2.6 we use the symbol $|\cdot|^{*},$ which means
$$\left|\sum_{M<m\leq 2M}z_m \right|^{*}=\max_{M<u\leq 2M} \left|\sum_{M<m\leq u}z_m\right|.$$

\section{\bf Some Lemmas}

In order to prove theorems, we need the following Lemmas.

{\bf Lemma 2.1.}  {\it For any $\mathcal{H}\geq 3$,
we have
\begin{eqnarray}
&&\psi(t)=-\sum_{1\leq |h|\leq \mathcal{H}}\frac{e(ht)}{2\pi hi}+O\left(\min\left(1,\frac{1}{\mathcal{H} \Vert t\Vert}\right)\right),\\
&&\min\left(1,\frac{1}{\mathcal{H} \Vert t\Vert}\right)=\sum_{h=-\infty}^\infty b(h)e(ht),
\end{eqnarray}
where
$$b(0)\ll \frac{\log \mathcal{H}}{\mathcal{H}}, \ \ b(h)\ll \min(1/|h|, \mathcal{H}/h^2).$$

}

\begin{proof}
 See Heath-Brown \cite{HB}.
\end{proof}

{\bf Lemma 2.2.} {\it Suppose $3< a<b\leq 2a$ and  the function $f(u)$ is at least $6$
times differentiable on the interval $[a,b]$ such that the estimate
$|f^{(j)}|\asymp Fa^{-j}\ (a\leq u\leq b)$
holds for some $F>0$.
  Then for any exponent pair $(\kappa, \lambda),$ we have the estimate
$$\sum_{a<n\leq b}e(f(n))\ll \frac aF+ F^\kappa a^{\lambda-\kappa}.$$
}

\begin{proof}
See, for example, Graham and Kolesnik \cite{GK}.
\end{proof}

{\bf Lemma 2.3.} {\it  Let $M>0,N>0,u_m>0,v_n>0,A_m>0,B_n>0(1\leq m\leq M,1\leq n\leq N),$ and let $Q_1$
and $Q_2$ be given non-negative numbers, $Q_1\leq Q_2$, then there is a $Q$
such that $Q_1\leq Q\leq Q_2$ and
\begin{eqnarray*}
  \sum_{m=1}^{M}A_{m}Q^{u_m}+\sum_{n=1}^{N}B_{n}Q^{-v_n}\ll \sum_{m=1}^{M}\sum_{n=1}^{N}
(A_{m}^{v_n}B_{n}^{u_m})^{\frac{1}{u_m+v_n}}\\
+\sum_{m=1}^{M}A_{m}Q_{1}^{u_m}+\sum_{n=1}^{N}B_{n}Q_{2}^{-v_n}.
\end{eqnarray*}

}

\begin{proof}
This is Lemma 2.4 of Graham and Kolesnik \cite{GK}.
\end{proof}

{\bf Lemma 2.4.} {\it Suppose $x^{1/2}<y<x^{2/3}$ is a parameter. Then
\begin{eqnarray}
E_f(x)=-\sum_{n\leq y}k_f(n)\psi\left(\frac xn\right)+ O(x^{1+\varepsilon}y^{-1}),
\end{eqnarray}
where $k_f(n): =f(n)-f(n-1)$ with $f(0)=0.$
}

\begin{proof}
This is contained in  Formula (3.7) of Zhai \cite{Zh2}.
\end{proof}

{\bf Lemma 2.5.}  {\it Suppose $M, N\geq 1$ are real numbers, $\alpha>0, \beta>0$ are fixed constants, $\Delta>0.$ Let
$\mathcal{G}(M,N;\Delta)$ denote the number of solutions of the inequality
$$\left|\left(\frac{n_1}{n_2}\right)^\alpha-\left(\frac{m_1}{m_2}\right)^\beta\right|\leq \Delta, \ \ n_1,n_2\sim N, m_1,m_2\sim M.$$
Then we have
$$\mathcal{G}(M,N;\Delta)\ll MN\log 2MN+\Delta M^2N^2,$$
where the $\ll$ constant is absolute.

}

\begin{proof}
This is Lemma 1 of Fouvry and Iwaniec \cite{FI}.
\end{proof}

{\bf Lemma 2.6.}  {\it  Let
$$S_1=\sum_{H<h\leq 2H}\sum_{N<n\leq 2N}a(h,n)\left|\sum_{M<m\leq 2M} e\left(U\frac{h^\beta n^\gamma m^\alpha}{H^\beta N^\gamma M^\alpha}\right)\right|^{*},$$
where $H, N, M$ are positive integers, $U$ is a real number greater than one, $a(h,n)$ is a complex number of modulus at most one; moreover,
$\alpha, \beta, \gamma$ are fixed real numbers such that $\alpha(\alpha-1)\beta\gamma\not= 0$. Then we have
\begin{equation}
S_1\ll (HNM)^{1+\varepsilon}\left(  \left(\frac{U}{HNM^2}\right)^{1/4} + \frac{1}{M^{1/2}} +  \frac{1}{U }    \right).
\end{equation}
}

\begin{proof}
 This is Theorem 3 of Robert and Sargos \cite{RS}.
\end{proof}

{\bf Lemma 2.7.}   {\it Let $\mathcal{H}\geq 3.$ Then $\psi(t)$ can be written as the form
$$\psi(t)=\sum_{1\leq |h|\leq \mathcal{H}}\alpha(h)e(ht)+O(\sum_{ |h|\leq \mathcal{H}}\beta(h)e(ht)),$$
where $\alpha(h)\ll 1/|h|,\ \beta(h)\ll 1/\mathcal{H}.$

}

\begin{proof}
See Vaaler \cite{V}, or  the appendix of Graham and Kolesnik \cite{GK}.
\end{proof}

\section{\bf Proof of Theorem 1}

Suppose $T\leq x\leq 2T$ and $T^{1/2}\ll y=y(T)\ll T^{2/3}$ is a parameter to be determined. By Lemma 2.4 we have
\begin{equation}
\int_T^{2T}|E_f(x)|^2dx\ll \int_1 +T^{3+\varepsilon}y^{-2},
\end{equation}
where
\begin{equation*}
 \int_1:=\int_T^{2T}\left|\sum_{n\leq y}k_f(n)\psi\left(\frac xn\right)\right|^2dx.
\end{equation*}

Take $\mathcal{H}=T^{2}$ in Lemma 2.1. We have
\begin{eqnarray}
&&\sum_{n\leq y}k_f(n)\psi\left(\frac xn\right)=\Sigma_1(x)+O(\Sigma_2(x)),\\
&&\Sigma_1(x)=-\sum_{1\leq |h|\leq \mathcal{H}}\frac{1}{2\pi ih}\sum_{n\leq y}k_f(n)e\left(\frac{hx}{n}\right),  \nonumber\\
&&\Sigma_2(x)=\sum_{n\leq y}|k_f(n)|\min\left(1,\frac{1}{\mathcal{H}\Vert \frac xn\Vert}\right).\nonumber
\end{eqnarray}

By Cauchy's inequality we have (note that $f(n)\ll n^\varepsilon$)
\begin{eqnarray}
\int_T^{2T}|\Sigma_2(x)|^2dx&&\ll \int_T^{2T}\sum_{n\leq y}|k_f(n)|^2\sum_{n\leq y}\min\left(1,\frac{1}{\mathcal{H}^2\Vert\frac xn\Vert^2}\right)dx\\
&&\ll y^{1+\varepsilon}\sum_{n\leq y} n\int_{T/n}^{2T/n}\min\left(1,\frac{1}{\mathcal{H}^2\Vert u\Vert^2}\right)du\nonumber\\
&&\ll Ty^{1+\varepsilon}\int_{-1/2}^{1/2}\min\left(1,\frac{1}{\mathcal{H}^2  u ^2}\right)du\nonumber\\
&&\ll Ty^{1+\varepsilon}\mathcal{H}^{-1}\ll 1.\nonumber
\end{eqnarray}

By a splitting argument we have
\begin{eqnarray}
\Sigma_1(x)\ll \frac{\log^2 T}{H}\left| \sum_{h\sim  H} \sum_{n\sim N}k_f(n)e\left(\frac{hx}{n}\right)\right|
\end{eqnarray}
for $1\ll H\ll \mathcal{H} $ and $1\ll N\ll y.$ So
\begin{eqnarray}
|\Sigma_1(x)|^2&&\ll \frac{\log ^4 T}{H^2}\sum_{ h_1,h_2\sim H } \sum_{n_1, n_2\sim  N}k_f(n_1)  k_f(n_2)
  e\left(x\left(\frac{h_1}{n_1}- \frac{h_2}{n_2}  \right)\right)\\
  &&=\frac{\log ^4 T}{H^2}\left(\Sigma_{11}(x)+\Sigma_{12}(x)\right),\nonumber
\end{eqnarray}
where
\begin{eqnarray*}
&&\Sigma_{11}(x):=\sum_{\stackrel{h_1,h_2\sim H;n_1, n_2\sim  N}{h_1n_2=h_2n_1}}k_f(n_1)  k_f(n_2) ,\\
&&\Sigma_{12}(x):=\sum_{\stackrel{h_1,h_2\sim H;n_1, n_2\sim  N}{h_1n_2\not=h_2n_1}}k_f(n_1)  k_f(n_2)
  e\left(x\left(\frac{h_1}{n_1}- \frac{h_2}{n_2}  \right)\right)
\end{eqnarray*}

Obviously we have
\begin{equation}
\int_T^{2T}\Sigma_{11}(x)dx\ll T^{1+\varepsilon}HN.
\end{equation}

By the first derivative test we have
\begin{eqnarray}
\int_T^{2T}\Sigma_{12}(x)dx&&\ll \sum_{\stackrel{h_1,h_2\sim H;n_1, n_2\sim  N}{h_1n_2\not=h_2n_1}}
\frac{ |k_f(n_1)  k_f(n_2)|  }{\left|\frac{h_1}{n_1}- \frac{h_2}{n_2}  \right|}\\
&&= \sum_{\stackrel{h_1,h_2\sim H;n_1, n_2\sim  N}{h_1n_2\not=h_2n_1}}
\frac{ |k_f(n_1)  k_f(n_2)n_1n_2|  }{\left|h_1n_2-h_2n_1 \right|}\nonumber\\
&&\ll H  N^{3+\varepsilon}.\nonumber
\end{eqnarray}

From (3.5)-(3.7) we get
\begin{eqnarray}
\int_T^{2T}|\Sigma_{1}(x)|^2dx\ll T^{1+\varepsilon} N+N^{3+\varepsilon}\ll T^{1+\varepsilon} y+y^{3+\varepsilon}.
\end{eqnarray}

Now Theorem 1 follows from (3.1)-(3.3)  and (3.8)  by choosing $y=T^{3/5}.$

\section{\bf Estimates of some exponential sums}

Let $x$ and $N$ be large real numbers with $x^{1/3}\ll N\ll x^{2/3}.$  Suppose $H\geq 1, M_1\geq 1, M_2\geq 2$ are real numbers.
Define
\begin{eqnarray*}
&&S_{\delta,I}(H, M_1,M_2):=\sum_{h\sim H}c_h\sum_{\stackrel{m_1\sim M_1,m_2\sim M_2}{m_1m_2\asymp N}}a_{m_1}e\left(\frac{hx}{m_1m_2+\delta}\right),\\
&&S_{\delta,II}(H, M_1,M_2):=\sum_{h\sim H}c_h\sum_{\stackrel{m_1\sim M_1,m_2\sim M_2}{m_1m_2\asymp N}}a_{m_1}b_{m_2}e\left(\frac{hx}{m_1m_2+\delta}\right),
\end{eqnarray*}
where $c_h, a_{m_1}, b_{m_2}\in {\Bbb C}$ such that $c_h\ll 1, a_{m_1}\ll m_1^{\varepsilon}, b_{m_2}\ll m_2^\varepsilon,$ and $0\leq \delta\leq 1.$

We first prove the following Lemma 4.1, which plays an important role in the proofs of both  Theorem 2 and Theorem 3.

{\bf Lemma 4.1.} {\it Suppose  $x^{1/3}\ll N\ll x^{2/3}, $ $N^{1/3}\ll M_1\ll N^{1/2}$ and $H\ll N^{1/2-\varepsilon}.$ For  any exponent pair $(\kappa, \lambda)$ we have the estimate
\begin{eqnarray}
S_{\delta,II}(H, M_1,M_2)N^{-\varepsilon} &&\ll H^{1/2}M_1^{1/2}M_2
 + Hx^{\frac{\kappa}{2}}M_1^{1-\kappa}M_2^{\frac{1+\lambda-2\kappa}{2}}\\
 &&\ \ \ +Hx^{\frac{\kappa}{2+2\kappa}}M_1^{\frac{2}{2+2\kappa}}M_2^{\frac{1+\lambda}{2+2\kappa}}+\frac{H^2x}{N^2}.\nonumber
\end{eqnarray}

Especially we have
 \begin{eqnarray}
  S_{\delta,II}(H, M_1,M_2)N^{-\varepsilon} &&\ll H^{1/2}M_1^{1/2}M_2
 + Hx^{\frac{1}{4}}M_1^{\frac{1}{2}}M_2^{\frac{1}{4}}\\&&
 +Hx^{\frac{1}{6}}M_1^{\frac{2}{3}}M_2^{\frac{1}{2}}+\frac{H^2x}{N^2}.\nonumber
\end{eqnarray}

}

\begin{proof}The idea of the proof of this lemma comes from \cite{HB}.
By Taylor's formula we have
$$
\frac{hx}{m_1m_2+\delta}=\frac{hx}{m_1m_2}-\frac{\delta hx}{m_1^2m_2^2}+O\left(\frac{hx}{M_1^3M_2^3}\right),
$$
which implies that
\begin{equation}
S_{\delta,II}(H, M_1,M_2)=\sum_{h\sim H}c_h\sum_{\stackrel{m_1\sim M_1,m_2\sim M_2}{m_1m_2  \asymp N}}
a_{m_1}b_{m_2}e\left(\frac{hx}{m_1m_2}- \frac{\delta hx}{m_1^2m_2^2}\right)+O\left(\frac{H^2x^{1+\varepsilon}}{N^2}\right).
\end{equation}
We only need to bound the sum
$$S_{\delta,II}^{*}(H, M_1,M_2):=\sum_{h\sim H}c_h\sum_{\stackrel{m_1\sim M_1,m_2\sim M_2}{m_1m_2 \asymp N}}
a_{m_1}b_{m_2}e\left(\frac{hx}{m_1m_2}- \frac{\delta hx}{m_1^2m_2^2}\right).$$

Obviously $0<h/m_1\leq 2H/M_1.$ Suppose $1\ll Q\ll HM_1$ is a positive integer parameter to be determined later.
For each $1\leq q\leq Q,$ let
$$E_q=\{(h, m_1):  h\sim H, m_1\sim M_1, \frac{2(q-1)H}{QM_1}<\frac{h}{m_1}\leq \frac{2qH}{QM_1}\}.$$

We can write $S_{\delta,II}^{*}(H, M_1,M_2)$ in the form
$$ S_{\delta,II}^{*}(H, M_1,M_2)=\sum_{m_2\sim M_2}b_{m_2}\sum_{q=1}^Q\sum_{\stackrel{(h,m_1)\in E_q}{m_1m_2\sim N}}c_h
a_{m_1} e\left(\frac{hx}{m_1m_2}- \frac{\delta hx}{m_1^2m_2^2}\right).$$
By Cauchy's inequality we get
\begin{eqnarray}
&&\ \ \ \ |S_{\delta,II}^{*}(H, M_1,M_2)|^2\\
&&\ll M_2^{1+\varepsilon}Q\sum_{m_2\sim M_2} \sum_{q=1}^Q \left|\sum_{\stackrel{(h,m_1)\in E_q}{m_1m_2 \asymp N}}c_h
a_{m_1} e\left(\frac{hx}{m_1m_2}- \frac{\delta hx}{m_1^2m_2^2}\right) \right|^2\nonumber\\
&&= M_2^{1+\varepsilon}Q\sum_{m_2\sim M_2} \sum_{q=1}^Q  \sum_{\stackrel{(h_1,m_{11})\in E_q}{m_{11}m_2 \asymp N}}
\sum_{\stackrel{(h_2,m_{12})\in E_q}{m_{12}m_2 \asymp N}} c_{h_1} a_{m_{11}}\overline{c_{h_2} a_{m_{12}}}e(r(m_2))\nonumber\\
&& =M_2^{1+\varepsilon}Q\sum_{q=1}^Q  \sum_{ (h_1,m_{11})\in E_q }
\sum_{ (h_2,m_{12})\in E_q } c_{h_1} a_{m_{11}}\overline{c_{h_2} a_{m_{12}}}
 \sum_{m_2\in I} e\left( r(m_2)\right),\nonumber
\end{eqnarray}
where $I=I(m_{11},m_{12})$ is a subinterval of $(M_2, 2M_2],$ and
\begin{eqnarray*}
r(m_2)=r(m_2;h_1,h_2,m_{11},m_{12}):=\frac{x}{m_2}\left(\frac{h_1}{m_{11}}-\frac{h_2}{m_{12}}\right)-
 \frac{\delta x}{m_2^2}\left( \frac{h_1}{m_{11}^2}-\frac{h_2}{m_{12}^2} \right).
 \end{eqnarray*}
Let
\begin{eqnarray*}
\eta=\frac{h_1}{m_{11}}-\frac{h_2}{m_{12}} .
 \end{eqnarray*}
It is easy to see that the contribution of $\eta=0$ to  $|S_{\delta,II}^{*}(H, M_1,M_2)|^2$
is $O(QM_1HM_2^2N^\varepsilon).$ Now we consider the case $\eta\not=0.$ From the conditions
$H\ll N^{1/2-\varepsilon}, N^{1/3}\ll M_1\ll N^{1/2}, M_1M_2\asymp N$, we see easily that
$$|r^{(j)}(u)|\asymp \frac{x |\eta|}{M_2^{j+1}},\ \ j=1, 2, 3, 4, 5, 6.$$

By Lemma 2.2 we see that the estimate
\begin{eqnarray}
\sum_{m_2\in I} e\left( r(m_2)\right)\ll \min\left(M_2, \frac{M_2^2}{x|\eta|}\right)+\left(\frac{x|\eta|}{M_2^2}\right)^{\kappa}M_2^\lambda
\end{eqnarray}
holds for any exponent pair $(\kappa,\lambda).$

By the definition of $E_q$ we have $|\eta|\leq 2H/QM_1.$ So from (4.4), (4.5) and the above discussions we get that
(we use a splitting argument to $1/|\eta|$)
\begin{eqnarray}
|S_{\delta,II}^{*}(H, M_1,M_2)|^2&&\ll QM_1HM_2^2N^\varepsilon+QM_2^{2+\varepsilon} \sum_{\stackrel{h_1,h_2\sim H; m_{11},m_{12}
\sim M_1}{|\eta|\leq M_2x^{-1}}} 1,\\
&&\ \ \ +QM_2^{3+\varepsilon}x^{-1}\sum_{\stackrel{h_1,h_2\sim H; m_{11},m_{12}
\sim M_1}{\frac{M_2}{x}\ll |\eta|\leq \frac{H}{QM_1}} }\frac{1}{|\eta|}\nonumber\\
&&\ \ \ +Q^{1-\kappa}x^{\kappa}H^{\kappa}M_1^{-\kappa}M_2^{1+\lambda-2\kappa} \sum_{\stackrel{h_1,h_2\sim H; m_{11},m_{12}
\sim M_1}{|\eta|\leq HQ^{-1}M_1^{-1}}} 1\nonumber\\
&&\ll QM_1HM_2^2N^\varepsilon+QM_2^{2+\varepsilon}  \mathcal{A}(H, M_1;M_2x^{-1}) \nonumber\\
&&\ \ \ +QM_2^{3+\varepsilon}x^{-1}\max_{\frac{M_2}{x}\ll   \Delta \leq \frac{H}{QM_1}}\frac 1\Delta\mathcal{A}(H,M_1;\Delta)
 \nonumber\\
&&\ \ \ +Q^{1-\kappa}x^{\kappa}H^{\kappa}M_1^{-\kappa}M_2^{1+\lambda-2\kappa}  \mathcal{A}(H,M_1;HQ^{-1}M_1^{-1})\nonumber
\end{eqnarray}
where $\mathcal{A}(H,M_1;\Delta)$ denotes the number of solutions of the inequality
$$|\eta|\leq \Delta, \ \ h_1,h_2\sim H, m_{11},m_{12}\sim M_1.$$
By Lemma 2.5 we get
\begin{equation}
\mathcal{A}(H,M_1;\Delta)\ll HM_1\log HM_1+\Delta HM_1^3.
\end{equation}
From (4.7) we have
\begin{eqnarray}
QM_2^{2+\varepsilon}  \mathcal{A}(H, M_1;M_2x^{-1})&&\ll (QHM_2^2M_1+QH(M_1M_2)^3x^{-1})N^{\varepsilon}\\
&& \ll QHM_2^2M_1N^{\varepsilon}\nonumber
\end{eqnarray}
and
\begin{eqnarray}
&&\ \ \ \ \ \ \ QM_2^{3+\varepsilon}x^{-1}\max_{\frac{M_2}{x}\ll | \Delta|\leq \frac{H}{QM_1}}\frac 1\Delta\mathcal{A}(H,M_1;\Delta)\\
&&\ll (QHM_2^2M_1+QH(M_1M_2)^3x^{-1})N^{\varepsilon}\ll  QHM_2^2M_1N^{\varepsilon}\nonumber
\end{eqnarray}
by noting that $N\ll x^{2/3},\ M_1\ll N^{1/2}, M_1M_2\asymp N.$ From (4.7) again we get
\begin{eqnarray}
&&\ \ \ \ Q^{1-\kappa}x^{\kappa}H^{\kappa}M_1^{-\kappa}M_2^{1+\lambda-2\kappa}  \mathcal{A}(H,M_1;HQ^{-1}M_1^{-1})\\
&&\ll Q^{1-\kappa}x^{\kappa}H^{\kappa}M_1^{-\kappa}M_2^{1+\lambda-2\kappa} (HM_1\log HM_1+Q^{-1}H^2M_1^2)\nonumber\\
&&\ll Q^{-\kappa}x^{\kappa}H^{2+\kappa}M_1^{2-\kappa}M_2^{1+\lambda-2\kappa}\nonumber
\end{eqnarray}
by recalling our assumption $Q\ll HM_1.$

From  (4.6) and  (4.8)-(4.10) we get
\begin{eqnarray*}
|S_{\delta,II}^{*}(H, M_1,M_2)|^2N^{-\varepsilon}&&\ll QHM_2^2M_1
 +Q^{-\kappa}x^{\kappa}H^{2+\kappa}M_1^{2-\kappa}M_2^{1+\lambda-2\kappa}.
\end{eqnarray*}
Choosing a best $Q\in [1, HM_1]$ via Lemma 2.3 we get
 \begin{eqnarray}
S_{\delta,II}^{*}(H, M_1,M_2)N^{-\varepsilon}&&\ll H^{1/2}M_1^{1/2}M_2
 + Hx^{\frac{\kappa}{2}}M_1^{1-\kappa}M_2^{\frac{1+\lambda-2\kappa}{2}}\\
 &&\ \ \ +Hx^{\frac{\kappa}{2+2\kappa}}M_1^{\frac{2}{2+2\kappa}}M_2^{\frac{1+\lambda}{2+2\kappa}}.\nonumber
\end{eqnarray}
  which combining  (4.3) gives (4.1). The estimate (4.2) follows from (4.1) by taking the exponent pair
  $(1/2,1/2).$
\end{proof}

{\bf Lemma 4.2.} {\it Suppose $x^{8/17}\ll N\ll x^{13/25}.$  If  $N^{1/3}\ll M_1\ll N^{1/2}$ and  $H\ll N^{1/2-\varepsilon},$ then
$$
S_{\delta,II}(H, M_1,M_2)  \ll H x^{\frac{47}{100}+\varepsilon}.
$$
}

\begin{proof} From (4.2) and  $N^{1/3}\ll M_1\ll N^{1/2}$ we see that
$$
S_{\delta,II}(H, M_1,M_2)N^{-\varepsilon} \ll H^{1/2}N^{5/6}
 + Hx^{\frac{1}{4}}N^{\frac{3}{8}}
 +Hx^{\frac{1}{6}}N^{\frac{7}{12}},
$$
which implies Lemma 4.2 by noting that $N\ll x^{13/25}.$
\end{proof}

{\bf Lemma 4.3.} {\it Suppose $x^{13/25}\ll N\ll x^{9/17}$. If $N^{2/9}\ll M_1\ll N^{4/9}$ and $H\ll N^{1/2-\varepsilon},$ then
$$
S_{\delta,II}(H, M_1,M_2)  \ll H x^{\frac{8}{17}+\varepsilon}.
$$
}

\begin{proof}
From (4.2) and  $N^{2/9}\ll M_1\ll N^{4/9}$ we see that
$$
S_{\delta,II}(H, M_1,M_2)N^{-\varepsilon} \ll H^{1/2}N^{8/9}
 + Hx^{\frac{1}{4}}N^{\frac{13}{36}}
 +Hx^{\frac{1}{6}}N^{\frac{31}{54}},
$$
which implies Lemma 4.3 by noting that $N\ll x^{9/17}.$
\end{proof}

{\bf Lemma 4.4.} {\it Suppose $x^{13/25}\ll N\ll x^{9/17}$. If $M_1\ll N^{5/9}$ and $H\ll Nx^{-8/17},$ then
$$
S_{\delta,I}(H, M_1,M_2)  \ll H x^{\frac{8}{17}+\varepsilon}.
$$
}

\begin{proof}
We consider three cases: $M_1\ll N^{2/9},$ $N^{2/9}\ll M_1\ll N^{4/9},$ $ N^{4/9}\ll M_1\ll N^{5/9}.$

If   $  M_1\ll N^{2/9}, $  then using the exponent pair $(2/7,4/7)$ to the sum over $m_2$ we get
\begin{eqnarray*}
S_{\delta,I}(H, M_1,M_2) && \ll H H^{2/7}x^{2/7+\varepsilon}M_1^{5/7}\ll H H^{2/7}x^{2/7+\varepsilon}N^{10/63}\\
&&\ll H (Nx^{-8/17})^{2/7}x^{2/7+\varepsilon}N^{10/63}\ll Hx^{18/119+\varepsilon}N^{4/9}\\
&&\ll Hx^{54/119+\varepsilon}\ll Hx^{8/17+\varepsilon}.
\end{eqnarray*}

When $N^{2/9} \ll M_1\ll N^{4/9},$ by Lemma 4.3 we get
 \begin{eqnarray*}
S_{\delta,I}(H, M_1,M_2)    \ll Hx^{8/17+\varepsilon}.
\end{eqnarray*}

 Finally suppose $ N^{4/9}\ll M_1\ll N^{5/9}.$
 Define
 $$g(n):=\sum_{\stackrel{n=m_1m_2}{m_1\sim M_1, m_2\sim M_2}}a_{m_1}.$$
 Then $S_{\delta,I}(H, M_1,M_2)$ can be written as
\begin{eqnarray}
S_{\delta,I}(H, M_1,M_2)&&=\sum_{h\sim H}c_h\sum_{M_1M_2<n\leq 4M_1M_2}g(n)e\left(\frac{hx}{n+\delta}\right)\\
&&=\sum_{h\sim H}c_h\sum_{M_1M_2<n\leq 4M_1M_2}g(n)e\left(\frac{hx}{n}\right) \Phi(n,h),\nonumber
\end{eqnarray}
 where
 $$\Phi(u,h):=e\left(\frac{hx}{u+\delta}-\frac{hx}{u}\right),\ \ M_1M_2<u\leq 4M_1M_2.$$

 It is easy to see that
 \begin{equation}
 \Phi(u,h)\ll 1,\ \ \frac{\partial}{\partial u}\Phi(u,h)\ll \frac{hx}{u^3}.
 \end{equation}

 Let
 $$A(u,h):=\sum_{M_1M_2<n\leq u}g(n)e\left(\frac{hx}{n}\right).$$
 By partial integration we have
 \begin{eqnarray}
 &&\sum_{M_1M_2<n\leq 4M_1M_2}g(n)e\left(\frac{hx}{n}\right) \Phi(n,h)
  =\int_{M_1M_2}^{4M_1M_2}\Phi(u,h)dA(u,h)\\
  &&=\Phi(4M_1M_2,h)A(4M_1M_2,h)-\int_{M_1M_2}^{4M_1M_2}A(u,h)\times \frac{\partial}{\partial u}\Phi(u,h)du\nonumber
 \end{eqnarray}

 From (4.12)-(4.14) we have
 \begin{eqnarray}
S_{\delta,I}(H, M_1,M_2)&&=\sum_{h\sim H}c_h \Phi(4M_1M_2,h)A(4M_1M_2,h) \\
&&\ \ \ \ -\int_{M_1M_2}^{4M_1M_2}\sum_{h\sim H}c_hA(u,h)\times \frac{\partial}{\partial u}\Phi(u,h)du\nonumber\\
&&\ll \left(1+\frac{Hx}{N^2}\right) \sum_{h\sim H}\sum_{m_1\sim M_1}a^{*}_{m_1}\left|\sum_{m_2\sim M_2}e\left(\frac{hx}{m_1m_2}\right)\right|^{*},\nonumber
\end{eqnarray}
 where $a^{*}_{m_1}=|a_{m_1}|.$ By Lemma 2.6 we have
 \begin{equation}
 x^{-\varepsilon}\sum_{h\sim H}\sum_{m_1\sim M_1}a^{*}_{m_1}\left|\sum_{m_2\sim M_2}e\left(\frac{hx}{m_1m_2}\right)\right|^{*}\ll
 Hx^{\frac 14}M_1^{\frac 12}M_2^{\frac 14}+HM_1M_2^{\frac 12}+\frac{N^2}{x}.
 \end{equation}

 From (4.15), (4.16) and the condition $N^{4/9}\ll M_2\ll N^{5/9}$ we get
 \begin{eqnarray}
x^{-\varepsilon}S_{\delta,I}(H, M_1,M_2)&& \ll  Hx^{\frac 14}M_1^{\frac 12}M_2^{\frac 14}+HM_1M_2^{\frac 12}+\frac{N^2}{x}\\
&&\ \ \ +\frac{H^2x^{\frac 54}}{N^2} M_1^{\frac 12}M_2^{\frac 14}+\frac{H^2x}{N^2}M_1M_2^{\frac 12}+H\nonumber\\
&&\ll Hx^{\frac 14}N^{\frac{7}{18}}+HN^{\frac 79}+\frac{N^2}{x}+
\frac{H^2x^{\frac 54}}{N^{\frac{29}{18}}}+\frac{H^2x}{N^{\frac{11}{9}}}\nonumber\\
&&\ll Hx^{\frac{7063}{15300}}\ll Hx^{\frac{8}{17}}\nonumber
\end{eqnarray}
 by noting that $x^{13/25}\ll N\ll x^{9/17}$ and $H\ll Nx^{-8/17}.$ This completes the proof of Lemma 4.4.
\end{proof}

Now we consider another exponential sum.
Let $x$ and $N$ be large real numbers with $x^{1/3}\ll N\ll x^{2/3}.$  Suppose $H\geq 1, D\geq 1, M_1\geq 1, M_2\geq 2$ are real numbers such that
$D^2M_1M_2\asymp N.$
Define
$$T_\delta(H, D, M_1,M_2):=\sum_{h\sim H}c_h\sum_{d\sim D}\rho_d\sum_{\stackrel{m_1\sim M_1,m_2\sim M_2}{d^2m_1m_2\asymp N}}a_{m_1}b_{m_2}e\left(\frac{hx}{d^2m_1m_2+\delta}\right),$$
where $c_h, \rho_d, a_{m_1}, b_{m_2}\in {\Bbb C}$ such that $c_h\ll 1, \rho_d\ll d^\varepsilon, a_{m_1}\ll m_1^{\varepsilon}, b_{m_2}\ll m_2^\varepsilon,$ and $0\leq \delta\leq 1.$ We will prove the following lemma  4.5,   which plays an important role in the proof of Theorem 3.

{\bf Lemma  4.5.} {\it Suppose $x^{1/3}\ll N\ll x^{2/3},$ $1\ll D\ll N^{1/6},$   $H \ll (N/D^2)^{1/2-\varepsilon}$ and    $(N/D^2)^{1/3}\ll M_1\ll (N/D^2)^{1/2}.$  Then for
any exponent pair $(\kappa, \lambda)$   we have
\begin{eqnarray}
   T_\delta(H, D, M_1,M_2)N^{-\varepsilon}&&\ll  \frac{H^{\frac 12}N^{\frac 56}}{D^{\frac 76}}+\frac{Hx^{\frac{\kappa}{2}}N^{\frac{3+
   \lambda-4\kappa}{4}}}{D^{\frac{1+\lambda-\kappa}{2}}}\\&&\ \ \ +
\frac{Hx^{\frac{\kappa}{2+2\kappa}}N^{\frac{3+
\lambda}{4+4\kappa}}}{D^{\frac{1+\lambda+\kappa }{2+2\kappa}}}+\frac{H^2x}{DN^2}. \nonumber
\end{eqnarray}

Especially we have
\begin{eqnarray}
   T_\delta(H, D, M_1,M_2)N^{-\varepsilon}&&\ll  \frac{H^{\frac 12}N^{\frac 56}}{D^{\frac 76}}+\frac{Hx^{\frac{1}{4}}N^{\frac{3 }{8}}}{D^{\frac{1 }{2}}}
   +\frac{Hx^{\frac{1}{6}}N^{\frac{7}{12}}}{D^{\frac{2 }{3}}}+\frac{H^2x}{DN^2}.
\end{eqnarray}

}

\begin{proof}
 Similar to (4.3) we have
\begin{equation}
T_\delta(H, D, M_1,M_2)=T_\delta^{*}(H, D, M_1,M_2) +O\left(\frac{H^2x^{1+\varepsilon}}{D^5M_1^2M_2^2}\right),
\end{equation}
where
$$T_\delta^{*}(H, D, M_1,M_2) :=\sum_{h\sim H}c_h\sum_{d\sim D}\rho_d\sum_{\stackrel{m_1\sim M_1,m_2\sim M_2}{m_1m_2 \asymp N}}
a_{m_1}b_{m_2}e\left(\frac{hx}{d^2m_1m_2}- \frac{\delta hx}{d^4m_1^2m_2^2}\right).$$

Obviously $0<h/d^2m_1\leq 2H/D^2M_1.$ Suppose $1\ll Q\ll DHM_1$ is a positive parameter to be determined later.
For each $1\leq q\leq Q,$ let
$$E_q^{*}=\{(h,d, m_1):  h\sim H,d\sim D, m_1\sim M_1, \frac{2(q-1)H}{QD^2M_1}<\frac{h}{d^2m_1}\leq \frac{2qH}{QD^2M_1}\}.$$

We now  write $T_\delta^{*}(H, D, M_1,M_2)$ in the form
$$ T_\delta^{*}(H,D, M_1,M_2)=\sum_{m_2\sim M_2}b_{m_2}\sum_{q=1}^Q\sum_{\stackrel{(h,d,m_1)\in E_q}{d^2m_1m_2\sim N}}c_h
\rho_d a_{m_1} e\left(\frac{hx}{d^2m_1m_2}- \frac{\delta hx}{d^4m_1^2m_2^2}\right).$$
Similar to (4.4)  we can  get
\begin{eqnarray}
\ \ \  |T_\delta^{*}(H, D, M_1, M_2)|^2\ll M_2^{1+\varepsilon}Q \sum_{\stackrel{h_j\sim H,d_j\sim D,m_{1j}\sim M_1(j=1,2)}{|\eta^{*}|\leq \frac{2H}{QD^2M_1}}}
\left| \sum_{m_2\in I^{*}} e\left( r_{*}(m_2)\right)\right|,
\end{eqnarray}
where $I^{*}=I^{*}(m_{11},m_{12})$ is a subinterval of $(M_2, 2M_2],$ and
\begin{eqnarray*}
r_*(m_2)&&=r(m_2;h_1,h_2,d_1,d_2,m_{11},m_{12}):=\frac{x}{m_2}\eta^{*}-
 \frac{\delta x}{m_2^2}\left( \frac{h_1}{d_1^4m_{11}^2}-\frac{h_2}{d_2^4m_{12}^2} \right),\\
\eta^{*}&&=\frac{h_1}{d_1^2m_{11}}-\frac{h_2}{d_2^2m_{12}} .
 \end{eqnarray*}

It is easy to see that the contribution of $\eta^{*}=0$ to  $|T_\delta^{*}(H, D, M_1, M_2)|^2$
is at most $O(QDM_1HM_2^2N^\varepsilon).$ Now we consider the case $\eta^{*}\not=0.$ From the conditions
$1\ll D\ll N^{1/6}, H\ll (N/D^2)^{1/2-\varepsilon},  D^2M_1M_2\asymp N$, we see easily that
$$|r_{*}^{(j)}(u)|\asymp \frac{x |\eta^{*}|}{M_2^{j+1}}\ (M_2\leq u\leq 2M_2),\ \ j=1, 2, 3, 4, 5, 6.$$
Suppose that $(\kappa,\lambda)$ is any exponent pair. By Lemma 2.2 we have the estimate
\begin{eqnarray}
\sum_{m_2\in I^{*}} e\left( r_*(m_2)\right)\ll \min\left(M_2, \frac{M_2^2}{x|\eta^{*}|}\right)+\left(\frac{x|\eta^{*}|}{M_2^2}\right)^{\kappa}M_2^\lambda.
\end{eqnarray}

By the definition of $E_q^{*}$ we have $|\eta^{*}|\leq 2H/QD^2M_1.$ So from (4.21), (4.22) and the above discussions we get that
(we use a splitting argument to $1/|\eta^{*}|$)
\begin{eqnarray}
&&\ \ \ \ \ \ \ \ \ \  |T_\delta^{*}(H, D, M_1,M_2)|^2\\&&\ll QM_1DHM_2^2N^\varepsilon+QM_2^{2+\varepsilon}  \mathcal{B}(H, D, M_1; M_2x^{-1}) \nonumber \\
&&\ \ \ +QM_2^{3+\varepsilon}x^{-1}\max_{\frac{M_2}{x}\ll   \Delta \leq \frac{2H}{QD^2M_1}}\frac 1\Delta\mathcal{B}(H, D, M_1; \Delta)
 \nonumber\\
&&\ \ \ +Q^{1-\kappa}x^{\kappa}H^{\kappa}D^{-2\kappa}M_1^{-\kappa}M_2^{1+\lambda-2\kappa}  \mathcal{B}(H, D, M_1; HQ^{-1}D^{-2}M_1^{-1})\nonumber
\end{eqnarray}
where $\mathcal{B}(H, D, M_1; \Delta)$ denotes the number of solutions of the inequality
\begin{equation}
\left|\frac{h_1}{d_1^2m_{11}}-\frac{h_2}{d_2^2m_{12}}\right|\leq \Delta, \ \ h_1,h_2\sim H, d_1,d_2\sim D,  m_{11},m_{12}\sim M_1.
\end{equation}

If (4.24) holds, then we have
\begin{eqnarray*}
\left|\frac{h_1m_{12}}{h_2m_{11}}-\frac{d_1^2}{d_2^2}\right|\leq \frac{8\Delta D^2M_1}{H},
\end{eqnarray*}
which combining Lemma 2.5 gives
\begin{eqnarray}
&&\ \ \ \ \ \mathcal{B}(H, D, M_1; \Delta)=\sum_{\stackrel{h_1,h_2\sim H, d_1,d_2\sim D,  m_{11},m_{12}\sim M_1}{\left|\frac{h_1}{d_1^2m_{11}}-\frac{h_2}{d_2^2m_{12}}\right|\leq \Delta}}1\\
&&\ll \sum_{\stackrel{h_1,h_2\sim H, d_1,d_2\sim D,  m_{11},m_{12}\sim M_1}{
 \left|\frac{h_1m_{12}}{h_2m_{11}}-\frac{d_1^2}{d_2^2}\right|\leq \frac{8\Delta D^2M_1}{H}}}1
 \ll  \sum_{\stackrel{n_1,n_2\sim HM_1, d_1,d_2\sim D}{
 \left|\frac{n_1}{n_2}-\frac{d_1^2}{d_2^2}\right|\leq \frac{8\Delta D^2M_1}{H}}}\tau(n_1)\tau(n_2)\nonumber\\
 &&\ll N^\varepsilon\sum_{\stackrel{n_1,n_2\sim HM_1, d_1,d_2\sim D}{
 \left|\frac{n_1}{n_2}-\frac{d_1^2}{d_2^2}\right|\leq \frac{8\Delta D^2M_1}{H}}}1\nonumber\\
 &&\ll(HDM_1+\Delta HM_1^3D^4)N^{\varepsilon},\nonumber
\end{eqnarray}
where we used the well-known bound $\tau(n)\ll n^\varepsilon.$

 From the conditions $1\ll D\ll N^{1/6}, x^{1/3}\ll N\ll x^{2/3},\   D^2M_1M_2\asymp N$ we get
 $$ (M_1M_2)^3D^4x^{-1}\ll M_2^2DM_1,$$
which combining (4.25) gives
\begin{eqnarray}
\ \ \ \ \ QM_2^{2+\varepsilon}  \mathcal{B}(H, D,  M_1; M_2x^{-1})&&\ll (QHM_2^2DM_1+QH(M_1M_2)^3D^4x^{-1})N^{\varepsilon}\\
&& \ll QHM_2^2DM_1N^{\varepsilon}\nonumber
\end{eqnarray}
and
\begin{eqnarray}
&&\ \ \ \ \ \ \ QM_2^{3+\varepsilon}x^{-1}\max_{\frac{M_2}{x}\ll | \Delta|\leq \frac{H}{QD^2M_1}}\frac 1\Delta\mathcal{B}(H, D, M_1; \Delta)\\
&&\ll (QHM_2^2DM_1+QH(M_1M_2)^3D^4x^{-1})N^{\varepsilon}\ll  QHM_2^2DM_1N^{\varepsilon}.\nonumber
\end{eqnarray}
  From (4.25) again we get
\begin{eqnarray}
&&\ \ \ \ Q^{1-\kappa}x^{\kappa}H^{\kappa}D^{-2\kappa}M_1^{-\kappa}M_2^{1+\lambda-2\kappa}  \mathcal{B}(H, D, M_1; HQ^{-1}D^{-2}M_1^{-1})\\
&&\ll Q^{1-\kappa}x^{\kappa}H^{\kappa}D^{-2\kappa}M_1^{-\kappa}M_2^{1+\lambda-2\kappa} (HM_1D +Q^{-1}H^2D^2M_1^2)N^\varepsilon\nonumber\\
&&\ll Q^{-\kappa}x^{\kappa}H^{2+\kappa}D^{2-2\kappa}M_1^{2-\kappa}M_2^{1+\lambda-2\kappa}\nonumber
\end{eqnarray}
by recalling our assumption $Q\ll HDM_1.$

From  (4.21) and  (4.26)-(4.28) we get
\begin{eqnarray*}
|T_\delta^{*}(H, D, M_1,M_2)|^2N^{-\varepsilon}&&\ll QHM_2^2DM_1
 +Q^{-\kappa}x^{\kappa}H^{2+\kappa}D^{2-2\kappa}M_1^{2-\kappa}M_2^{1+\lambda-2\kappa}.
\end{eqnarray*}
Choosing a best $Q\in [1, HDM_1]$ via Lemma 2.3 we get
 \begin{eqnarray}
 && \ \ \ \ \ \  T_\delta^{*}(H, D, M_1,M_2)N^{-\varepsilon} \\
 &&\ll H^{1/2}M_2M_1^{1/2}D^{1/2}
 + Hx^{\frac{\kappa}{2}}D^{1-\frac{3\kappa}{2}}M_1^{1-\kappa}M_2^{\frac{1+\lambda-2\kappa}{2}} \nonumber\\
 &&\ \ \ +Hx^{\frac{\kappa}{2+2\kappa}}D^{\frac{2-\kappa}{2+2\kappa}}M_1^{\frac{2}{2+2\kappa}}M_2^{\frac{1+\lambda}{2+2\kappa}}\nonumber\\
 &&\ll  \frac{H^{\frac 12}N^{\frac 56}}{D^{\frac 76}}+\frac{Hx^{\frac{\kappa}{2}}N^{\frac{3+\lambda-4\kappa}{4}}}{D^{\frac{1+\lambda-\kappa}{2}}}+
\frac{Hx^{\frac{\kappa}{2+2\kappa}}N^{\frac{3+\lambda}{4+4\kappa}}}{D^{\frac{1+\lambda+\kappa }{2+2\kappa}}}\nonumber
\end{eqnarray}
by recalling that $(N/D^2)^{1/3}\ll M_1\ll (N/D^2)^{1/2}$ and $D^2M_1M_2\asymp N.$
Now (4.18) follows from (4.20) and (4.29). The estimate (4.19) follows from (4.18) by taking the exponent pair $(1/2,1/2)$.
\end{proof}

{\bf Remark. } When we use Lemma 4.1 and Lemma 4.5 to prove Theorem 2 and Theorem 3, we choose the exponent pair $(1/2,1/2).$
If we choose exponent pairs as well as we can, then we can slightly improve both the exponent  $8/17$ in Theorem 2 and the exponent
$73/155$ in Theorem 3. Especially, if $(\varepsilon, 1/2+\varepsilon)$ is an exponent pair, then the exponent  $8/17$ in Theorem 2 can be
improved to $7/15.$

\section{\bf Proof of Theorem 2}

\subsection{\bf Proof of Theorem 2}

Suppose $k\geq 2$ and $\ell\geq 1$ are fixed integers and  $f\in \{\tau_k, \Lambda_\ell, \mu_\ell, \omega_\ell\}.$
Let $  y= x^{9/17}, $ then $x/y=x^{8/17}.$    By Lemma 2.4 we have
\begin{eqnarray}
\ \ \ \ \ \ E_f(x)&&=-\sum_{n\leq  x^{9/17}}k_f(n)\psi\left(\frac xn\right)+ O(x^{8/17+\varepsilon})\\
&& =-\sum_{x^{8/17}<n\leq x^{9/17}}f(n)\psi\left(\frac xn\right)+ \sum_{x^{8/17}<n\leq x^{9/17}}f(n)\psi\left(\frac{x}{n+1}\right)\nonumber\\&&\ \ \ \ \ \
 + O(x^{8/17+\varepsilon} ).\nonumber
\end{eqnarray}
So we only need to bound the sum
\begin{eqnarray*}
R_{f,\delta}(N;x):=\sum_{n\sim N}f(n) \psi\left(\frac{x}{n+\delta}\right)\ \ (\delta=0,1)
\end{eqnarray*}
for  $x^{8/17}\ll N\ll x^{9/17}.$  Let $  \mathcal{H}: = Nx^{-8/17}.$   By (2.1) of Lemma 2.1 we have
\begin{equation}
R_{f,\delta}(N;x)=\Sigma_{f,\delta 1}+O(\Sigma_{f,\delta 2}),
\end{equation}
where
\begin{eqnarray*}
&&\Sigma_{f,\delta 1}=-\sum_{1\leq |h|\leq \mathcal{H}}\frac{1}{2\pi ih}\sum_{n\sim N}f(n)e\left(\frac{hx}{n+\delta}\right),\\
&&\Sigma_{f,\delta 2}=\sum_{n\sim N}|f(n)|\min\left(1,\frac{1}{ \mathcal{H}\Vert \frac{x}{n+\delta}\Vert}\right).
\end{eqnarray*}

By (2.2) of  Lemma 2.1 and Lemma 2.2 with the exponent pair $(2/7,4/7)$ we have
\begin{eqnarray}
\Sigma_{f,\delta 2}&&\ll N^\varepsilon \sum_{n\sim N} \min\left(1,\frac{1}{ \mathcal{H}\Vert \frac{x}{n+\delta}\Vert}\right)\\
&&\ll N^\varepsilon\left( \frac{N}{\mathcal{H}}+\sum_{h=1}^\infty  \min(h^{-1}, \mathcal{H}h^{-2})\sum_{n\sim N}e\left(\frac{hx}{n+\delta}\right)\right)\nonumber\\
&&\ll N^\varepsilon\left( \frac{N}{\mathcal{H}}+ \mathcal{H}^{2/7}x^{2/7}+\frac{N^2}{x}\right)\ll x^{8/17+\varepsilon}.\nonumber
\end{eqnarray}

By a splitting argument we get
\begin{eqnarray}
&&\Sigma_{f,\delta 1}\ll \frac 1H\left|\sum_{h\sim H} \sum_{n\sim N}f(n)e\left(\frac{hx}{n+\delta}\right)\right|\log x
\end{eqnarray}
for some $1\ll H\ll \mathcal{H}.$     So Theorem 2 follows from (5.1)-(5.4)  and the estimate
\begin{eqnarray}
S_{f,\delta}(H, N):= \sum_{h\sim H} \sum_{n\sim N}f(n)e\left(\frac{hx}{n+\delta}\right)\ll Hx^{8/17+\varepsilon}.
\end{eqnarray}

We only need to prove that (5.5) holds for $f\in \{\tau_k, \Lambda_\ell, \mu_\ell, \omega_\ell\}.$
We can prove the following two lemmas.

{\bf Lemma 5.1.}  {\it  The sum $S_{f,\delta}(H, N)$
can be written as  a sum of   $(\log^{\nu_f} x)$ expressions  of the form(Type I sum)
\begin{equation}
S_{f;I}(H,M_1,M_2):=\sum_{h\sim H}c_h\sum_{\stackrel{m_1\sim M_1, m_2\sim M_2}{m_1m_2 \asymp N}}a_{m_1}e\left(\frac{hx}{m_1m_2+\delta}\right)
\end{equation}
with $M_1\ll N^{1/3}$ and  expressions of the form(Type II sum)
\begin{equation}
S_{f;II}(H,M_1,M_2):=\sum_{h\sim H}c_h\sum_{\stackrel{m_1\sim M_1, m_2\sim M_2}{m_1m_2 \asymp N}}a_{m_1}b_{m_2}e\left(\frac{hx}{m_1m_2+\delta}\right)
\end{equation}
with $N^{1/3}\ll M_1\ll N^{1/2}, $ where $\nu_f\geq 1$ is a fixed integer, and $c_h\ll 1, a_{m_1}\ll m_1^{\varepsilon}, b_{m_2}\ll m_2^{\varepsilon}.$
Note that when $f=\omega_\ell$, an additional term $HN^{2/3}$ should be added to the above sums. }

{\bf Lemma 5.2.}  {\it  The sum $S_{f,\delta}(H, N)$
can be written as a sum of $(\log^{\nu_f} x)$ expressions of the form (5.6)
with $M_1\ll N^{5/9}$ and expressions of the form (5.7)
with $N^{2/9}\ll M_1\ll N^{4/9}, $ where $ \nu_f\geq 1$ is a fixed integer, and $c_h\ll 1, a_{m_1}\ll m_1^{\varepsilon}, b_{m_2}\ll m_2^{\varepsilon}.$ Note that when $f=\omega_\ell$, an additional term $HN^{2/3}$ should be added to the above sums. }

\bigskip

First suppose $x^{8/17}\ll N\ll x^{13/25}.$  By Lemma 4.2 we have
\begin{equation}
S_{f;II}(H,M_1,M_2)\ll Hx^{\frac{47}{100}+\varepsilon}.
\end{equation}

Using Lemma 2.2 with the exponent pair $(2/7,4/7)$ to estimate the sum over $m_2$ and estimate the sums over $h$ and $m_1$ trivially we get
\begin{eqnarray}
x^{-\varepsilon}S_{f;I}(H,M_1,M_2)&&\ll HH^{\frac 27}x^{\frac 27}M_1^{\frac 57}+\frac{N^2}{x}\\
&&\ll H(Nx^{-8/17})^{\frac 27}x^{\frac 27}N^{\frac{5}{21}}+\frac{N^2}{x}\nonumber\\&&\ll Hx^{3/7}\ll Hx^{8/17}\nonumber
\end{eqnarray}
by noting that $H\ll Nx^{-8/17}$ and $M_1\ll N^{1/3}.$  So from (5.8)  and (5.9) we get an estimate better than (5.5) in the  range
$x^{8/17}\ll N\ll x^{13/25}.$

When $x^{13/25}\ll N\ll x^{9/17},$ the estimate (5.5) follows  from Lemma 4.3, Lemma 4.4 and Lemma 5.2. This completes the proof of Theorem 2.

We will prove Lemma 5.1 and Lemma 5.2 in the next two subsections.

\subsection{\bf Proof of Lemma 5.1.}

In this subsection we will prove Lemma 5.1. We consider only $f\in \{\tau_k, \Lambda_\ell,\omega_\ell\},$
since the proofs for $\Lambda_\ell$ and $\mu_\ell$ are   the same.

We first prove the following decomposition formula.

{\bf Lemma 5.3.}  {\it Suppose $K\geq 2$ is a fixed integer,  $N_1\geq 1,\cdots, N_K\geq 1$ are natural numbers such that
$N_1\cdots N_K\asymp N$,  $W(u)$ is any function defined on $(N, 2N],$  $a_j(n_j)\in {\Bbb C} $
such that $a_j(n_j)\ll n_j^\varepsilon\ (j=1,\cdots, K).$ If $\max(N_1,\cdots, N_K)\ll N^{2/3},$ then the sum
\begin{equation}
\sum_{n_1\sim N_1}a_1(n_1)\cdots \sum_{n_K\sim N_K}a_K(n_K)W(n_1\cdots n_K)
\end{equation}
can be written as the form
\begin{equation}
\sum_{m_1\sim M_1}a_{m_1}\sum_{m_2\sim M_2}b_{m_2}W(m_1m_2)
\end{equation}
such that $N^{1/3}\ll M_1\ll N^{1/2}$ and $a_{m_1}\ll m_1^\varepsilon, \ b_{m_2}\ll m_2^\varepsilon. $
}

\begin{proof} Without loss of generality, suppose $N_1\leq N_2\leq \cdots \leq N_K.$

If $K=2,$ then $N_1\leq N_2\ll N^{2/3}$ and $N_1N_2\asymp N$ implies that $N^{1/3}\ll N_1\ll N^{1/2}.$ So we see that
(5.10) is of the form (5.11) by taking $m_1=n_1, m_2=n_2.$

Now suppose $K\geq 3$. We consider three cases:
$N_K<N^{1/3}$, $N^{1/3}\leq N_K\ll N^{1/2},$ $N^{1/2}\ll N_K\ll N^{2/3}.$

{\bf Case I.}  $N_K<N^{1/3},$ which implies that  $N_j<N^{1/3}(j=1,\cdots, K).$    Since $\prod_{j=1}^K N_j\asymp N, $
we find that there is a $j$ such that $2\leq j<K$ and
$$N_1\cdots N_{j-1}<N^{1/3},\ \ N_1\cdots N_j\gg N^{1/3}.$$
Thus we have
$$N^{1/3}\ll  N_1\cdots N_j=N_1\cdots N_{j-1}\times N_j \ll N^{2/3}.$$
If $N_1\cdots N_j\ll N^{1/2},$ then  the sum   (5.10) can be written as (5.11) by taking
\begin{eqnarray*}
&&m_1=n_1\cdots n_j ,\ \ m_2=n_{j+1}\cdots n_{K},\ M_1=N_1\cdots N_j ,\ \ M_2=N_{j+1}\cdots N_{K}  ,  \\
&& a_{m_1}= a_1(n_1)\cdots a_j(n_j)\ll m_1^\varepsilon,\ \ b_{m_2}=a_{j+1}(n_{j+1})\cdots a_{K}(n_{K})\ll m_2^\varepsilon.
 \end{eqnarray*}
 If $N_1\cdots N_j\gg N^{1/2},$ then  the sum   (5.10) can be written as (5.11) by taking
\begin{eqnarray*}
 m_1=  n_{j+1}\cdots n_{K}, \ \ m_2=n_1\cdots n_j .
 \end{eqnarray*}

{\bf Case II.}  $N^{1/3}\ll N_K\ll N^{1/2}. $

In this case  the sum (5.10)   can be written as (5.11) by taking
\begin{eqnarray*}
&&m_1=n_K ,\ \ m_2=n_1\cdots n_{K-1},\  M_1=N_K, \, M_2=N_1\cdots N_{K-1},  \\
&& a_{m_1}=a_K(n_K)\ll m_1^\varepsilon,\ \ b_{m_2}=a_1(n_1)\cdots a_{K-1}(n_{K-1})\ll m_2^\varepsilon.
 \end{eqnarray*}

{\bf Case III.}  $N^{1/2}\ll N_K\ll N^{2/3} . $

In this case  the sum   (5.10) can be written as (5.11) by taking
\begin{eqnarray*}
&&m_1=n_1\cdots n_{K-1} ,\ \ m_2=  n_{K},\ M_1=N_1\cdots N_{K-1} ,\ \ M_2= N_{K}  ,  \\
&& a_{m_1}= a_1(n_1)\cdots a_{K-1}(n_{K-1})\ll m_1^\varepsilon,\ \ b_{m_2}=  a_{K}(n_{K})\ll m_2^\varepsilon.
 \end{eqnarray*}

This completes the proof of Lemma 5.3.
\end{proof}

\subsubsection{Proof of Lemma 5.1 for $\tau_k$}

 Suppose $k\geq 2$ is a fixed integer. By a splitting argument we see that
$S_{\tau_k,\delta}(H, N)$ can be written as a sum of $O(\log^{k-1} N)$ exponential sums of the form
 \begin{eqnarray*}
S_{\tau_k,\delta}(H, N_1,\cdots, N_k):= \sum_{h\sim H} \sum_{\stackrel{n_1\sim N_1,\cdots, n_k\sim N_k}{n_1\cdots n_k \asymp N}}
e\left(\frac{hx}{n_1\cdots n_k+\delta}\right),
\end{eqnarray*}
where $N_1\geq 1, \cdots, N_k\geq 1$ are natural numbers such that $N_1  \cdots N_k\asymp N.$ Without loss of generality, we suppose
$N_1\leq N_2\leq \cdots \leq N_k.$

If  $N_k\ll N^{2/3},$ then from Lemma 5.3 we see that $S_{\tau_k,\delta}(H, N_1,\cdots, N_k)$ can be written as the form
(5.7) with $N^{1/3}\ll M_1\ll N^{1/2}.$ If $N_k\gg N^{2/3},$ then $N_1\cdots N_{k-1}\ll N^{1/3}.$ So $S_{\tau_k,\delta}(H, N_1,\cdots, N_k)$ can be written as the form
(5.6) with $  M_1\ll N^{1/3}.$

 \subsubsection{Proof of Lemma 5.1 for $f=\Lambda_\ell, \ell=1$}
 We first consider the function $\Lambda_\ell $ for the case $\ell=1$.
 Let  $u\geq 3$ be a fixed integer. If $2v^u\geq n$, then we have
 Heath-Brown's identity
 $$\Lambda(n)=\sum_{j=1}^u (-1)^j\left(_j^u\right)\sum_{\stackrel{n=\prod_{t=1}^jn_t\prod_{t=1}^jn_{j+t}}{n_1,\cdots, n_j<v}}\mu(n_1)\cdots \mu(n_j)\log n_{2u}.$$
 Using Heath-Brown's identity with $u=4$  we find that
 $S_{\Lambda,\delta}(H, N)$ can be written as a sum of $O(\log^7 N)$ exponential sums of the form
 \begin{eqnarray*}
S_{\Lambda,\delta}(H, N_1,\cdots, N_8):= \sum_{h\sim H} \sum_{\stackrel{n_1\sim N_1,\cdots, n_8\sim N_8}{n_1\cdots n_8\asymp N}}\log n_8  \prod_{j=1}^4\mu(n_j)
e\left(\frac{hx}{n_1\cdots n_8+\delta}\right),
\end{eqnarray*}
where $N_1\geq 1, \cdots, N_8\geq 1$ are natural numbers such that $N_1  \cdots N_8\asymp N, \ N_j\leq N^{1/4}(j=1,2,3,4).$

 If $\max(N_1,\cdots, N_8)\ll N^{2/3},$ then from Lemma 5.3 we see that $S_{\Lambda,\delta}(H, N_1,\cdots, N_8)$ can be written as the form (5.7)
 with  $N^{1/3}\ll M_1\ll N^{1/2}.$ Suppose now $N_j=\max(N_1,\cdots, N_8)\gg N^{2/3},$  then it is easy to see that $j\in \{5,6,7,8\}.$ So the sum  $S_{\Lambda,\delta}(H, N_1,\cdots, N_8)$ can be written as the form (5.6) with $  M_1\ll N^{1/3}.$

 \subsubsection{Proof of Lemma 5.1 for  $\Lambda_\ell, \ell\geq 2$}

 Since $\Lambda_\ell(n)=\sum_{n=n_1\cdots n_\ell}\Lambda(n_1)\cdots \Lambda(n_\ell),$ the sum
  $S_{\Lambda_\ell,\delta}(H, N)$ can be written as a sum of $O(\log^\ell N)$ exponential sums of the form
 \begin{eqnarray*}
S_{\Lambda_\ell,\delta}(H, N_1,\cdots, N_\ell):= \sum_{h\sim H} \sum_{\stackrel{n_1\sim N_1,\cdots, n_\ell\sim N_\ell}{n_1\cdots n_\ell\asymp N}} \prod_{j=1}^\ell\Lambda(n_j)
e\left(\frac{hx}{n_1\cdots n_\ell+\delta}\right),
\end{eqnarray*}
where $N_1\geq 1, \cdots, N_\ell\geq 1$ are natural numbers such that $N_1  \cdots N_\ell\asymp N.$ Without loss of generality, suppose
$N_1\leq N_2\cdots \leq N_\ell.$

If   $N_\ell \leq N^{2/3},$ then from Lemma 5.3 we see that   $S_{\Lambda_\ell,\delta}(H, N_1,\cdots, N_\ell)$ can be written as the form (5.7)
 with  $N^{1/3}\ll M_1\ll N^{1/2}.$

Now suppose $N_\ell>N^{2/3}.$ From $N_1\cdots N_\ell\asymp N,$ we get that $N_1  N_2\cdots  N_{\ell-1}\ll N^{1/3}.$
Applying  Heath-Brown's identity with $u=4$ to $\Lambda(n_\ell)$ again we find that $S_{\Lambda_\ell,\delta}(H, N_1,\cdots, N_\ell) $
can be written as a sum of $O(\log^7 N)$ exponential sums of the form
 \begin{eqnarray*}
 S_{\Lambda_\ell,\delta}(H, {\bf N}):&&= \sum_{h\sim H}
 \sum_{ \stackrel{n_j\sim N_j}{1\leq j\leq \ell-1} } \prod_{j=1}^{\ell-1}\Lambda(n_j)\\
 &&\ \ \
 \sum_{\stackrel{n_{\ell 1}\sim N_{\ell 1},\cdots, n_{\ell 8}\sim N_{\ell 8}}{n_{\ell 1}\cdots n_{\ell 8}\asymp N_\ell}} \prod_{j=1}^4\mu(n_{\ell j})
 \log n_{\ell 8}\times
e\left(\frac{hx}{n_1\cdots n_{\ell-1}n_{\ell 1}\cdots n_{\ell 8}+\delta}\right),
\end{eqnarray*}
 where ${\bf N}=(N_1, \cdots, N_{\ell-1},N_{\ell 1},\cdots, N_{\ell 8}),$ $N_{\ell 1}\geq 1,\cdots, N_{\ell 8}\geq 1$ are natural numbers
 such that $N_{\ell 1} \cdots N_{\ell 8}\asymp N_\ell,$ and $N_{\ell j}\leq N_{\ell}^{1/4}\ (j=1,2,3,4).$

 If there exists an $N_{\ell j}$ such that $N_{\ell j}\gg N^{2/3},$ then we have $5\leq j\leq 8.$ So the sum
 $S_{\Lambda_\ell,\delta}(H, {\bf N})$ can be written as the form (5.6) with $M_1\ll N^{1/3}.$ If $N_{\ell j}\ll N^{2/3}(1\leq j\leq 8),$ then
 by Lemma 5.3 the sum $S_{\Lambda_\ell,\delta}(H, {\bf N})$ can be written as the form (5.7) with $N^{1/3}\ll M_1\ll N^{1/2}.$

 \subsubsection{ Proof of Lemma 5.1 for $f=\omega_\ell, \ell=1$}

  Since $\omega(n)=\sum_{n=pn_1}1,$ the sum
 $S_{\omega,\delta}(H, N)$ can be written as a sum of $O(\log N)$ exponential sums of the form
 \begin{eqnarray*}
S_{\omega,\delta}(H, P, N_1):= \sum_{h\sim H} \sum_{\stackrel{n_1\sim N_1,  p\sim P}{n_1p\asymp N}}
e\left(\frac{hx}{n_1p+\delta}\right),
\end{eqnarray*}
where $N_1\geq 1,   P\geq 1$ are natural numbers such that $N_1 P\asymp N.$

If $N_1\gg N^{2/3},$  then the sum $S_{\omega,\delta}(H, P, N_1)$ can be written as the form (5.6) by taking
$m_1=p,m_2=n_1.$ If $N^{1/2}\ll N_1\ll N^{2/3},$  then  $S_{\omega,\delta}(H, P, N_1)$ can be written as the form (5.7) by taking
$m_1=p,m_2=n_1.$ If $N^{1/3}< N_1\ll N^{1/2},$  then   $S_{\omega,\delta}(H, P, N_1)$ can be written as the form (5.7) by taking
$m_1=n_1,m_2=p.$

Now suppose $N_1\ll N^{1/3},$ namely, $P\gg N^{2/3}.$ Then we have
\begin{eqnarray}
S_{\omega,\delta}(H, P, N_1)\ll \frac{1}{\log P}\left|S_{\omega,\delta}^{*}(H, P, N_1)\right|+HN_1P^{1/2},
\end{eqnarray}
 where
  \begin{eqnarray*}
S_{\omega,\delta}^{*}(H, P, N_1):= \sum_{h\sim H} \sum_{\stackrel{n_1\sim N_1,  n\sim P}{n_1n\asymp N}}\Lambda(n)
e\left(\frac{hx}{n_1n+\delta}\right).
\end{eqnarray*}

 We use Heath-Brown's identity with $u=4$ to $\Lambda(n).$ So we find that
 $S_{\omega,\delta}^{*}(H, P, N_1) $ can be written as a sum of $O(\log^7 N)$ exponential sums of the form
 \begin{eqnarray*}
S_{\omega,\delta}^{*}(H, N_1,N_2, \cdots, N_9):= \sum_{h\sim H} \sum_{\stackrel{n_1\sim N_1,\cdots, n_9\sim N_9}{n_1\cdots n_9 \asymp N}}\log n_9 \prod_{j=2}^5\mu(n_j)
e\left(\frac{hx}{n_1\cdots n_9+\delta}\right),
\end{eqnarray*}
where $N_2\geq 1, \cdots, N_9\geq 1$ are natural numbers such that $N_2  \cdots N_9\asymp P, \ N_j\leq P^{1/4}(j=2,3,4,5).$

If there exists an $N_j$ for which $N_j\gg N^{2/3},$ then we have $6\leq j\leq 9.$
Without loss of generality, suppose $j=9.$ We see that
  the sum $S_{\omega,\delta}^{*}(H, N_1,N_2, \cdots, N_9) $ can be written as the form (5.6) with $m_1\ll N^{1/3}$ by taking
$m_1=n_1\cdots n_8,m_2=n_9.$ If $N_j\ll N^{2/3}(j=2\leq j\leq 9),$ then from Lemma 5.3 we see that
 $S_{\omega,\delta}^{*}(H, N_1,N_2, \cdots, N_9) $ can be written as the form (5.7) with $N^{1/3}\ll m_1\ll N^{1/2}.$

Note that $HN_1P^{1/2}$ in (5.12) satisfies
$HN_1P^{1/2}\ll HN^{2/3}.$

\subsubsection{Proof of Lemma 5.1 for  $f=\omega_\ell, \ell\geq 2$}

  Since $\omega_\ell(n)=\sum_{n=n_1\cdots n_\ell}\omega(n_1)\cdots \omega(n_\ell),$ the sum
  $S_{\omega_\ell,\delta}(H, N)$ can be written as a sum of $O(\log^\ell N)$ exponential sums of the form
 \begin{eqnarray*}
S_{\omega_\ell,\delta}(H, N_1,\cdots, N_\ell):= \sum_{h\sim H} \sum_{\stackrel{n_1\sim N_1,\cdots, n_\ell\sim N_\ell}{n_1\cdots n_\ell \asymp N}} \prod_{j=1}^\ell\omega(n_j)
e\left(\frac{hx}{n_1\cdots n_\ell+\delta}\right),
\end{eqnarray*}
where $N_1\geq 1, \cdots, N_\ell\geq 1$ are natural numbers such that $N_1  \cdots N_\ell\asymp N$ and
$N_1\leq N_2\cdots \leq N_\ell.$

If $N_\ell\leq N^{2/3},$  then   by Lemma 5.3 we see that
the sum $S_{\omega_\ell,\delta}(H, N_1,\cdots, N_\ell)$ can be written as the form (5.7) with $N^{1/3}\ll M_1\ll N^{1/2}.$

If $N_\ell\gg N^{2/3},$ then the sum $\sum_{n_\ell\sim N_\ell}$ can be written as sums of the form
$\sum_{p\sim P_\ell}\sum_{n\sim N_\ell^{*}}$ with $P_\ell N_\ell^{*}\asymp N_\ell.$
If $N_\ell^{*}\gg N^{2/3},$ we get a sum of the form (5.6) with $M_1\ll N^{1/3}.$ If
$N^{1/3}\ll N_\ell^{*}\ll N^{2/3},$ we can  get a sum of the form (5.7) with $N^{1/3}\ll M_1\ll N^{1/2}.$
Finally suppose  $N_\ell^{*}\ll N^{1/3}. $     If $P_\ell \ll N^{2/3},$ then we use Lemma 5.3 directly to get a sum (5.7) with
$N^{1/3}\ll M_1\ll N^{1/2}.$
If $P_\ell \gg N^{2/3},$ then we apply Henth-Brown's identity to $p$ and then use  the
 procedure of the proof of $\omega_1.$ We omit the details.

\subsection{\bf Proof of Lemma 5.2}

In this subsection we will prove Lemma 5.2. Since the proof is similar to that of Lemma 5.1, we give only
a short description.

We first give  the following decomposition formula without proof.

{\bf Lemma 5.4.}  {\it Under the conditions of Lemma 5.3, if $\max(N_1,\cdots, N_K)\ll N^{4/9},$ then the sum (5.10) can
  be written as the form (5.11)
such that $N^{2/9}\ll M_1\ll N^{4/9}$ and $a_{m_1}\ll m_1^\varepsilon, \ b_{m_2}\ll m_2^\varepsilon. $
}

\bigskip

From the proof of Lemma 5.1 we see that if $f\in \{\tau_k, \Lambda_\ell, \mu_\ell, \omega_\ell\},$ then the sum
$S_{f,\delta}(H,N)$ can be written as a sum of $O((\log x)^{\nu_f})$ expressions of the form
$$S_{f;\delta}(H,N_1,\cdots, N_K):=\sum_{h\sim H}c_h\sum_{n_1\sim N_1}a_{n_1}\cdots \sum _{n_K\sim N_K}a_{n_K}e\left(\frac{hx}{n_1\cdots n_K+\delta}\right),$$
where $\nu_f\geq 1$ and $K=K(f)\geq 2$ are fixed integers, $N_1\geq 1,\cdots, N_K\geq 1$ and $N_1\cdots N_K\asymp N.$

Let $N_j=\max(N_1,\cdots, N_K).$
If   $N_j\gg N^{4/9},$ then $a_{n_j}=1$ or $a_{n_j}=\log n_j.$ Hence  the sum $S_{f;\delta}(H,N_1,\cdots, N_K)$ can be written as the form
(5.6) with $m_1\ll N^{5/9}.$  If $N_j\ll N^{4/9}$, then from Lemma 5.4  we see that the sum $S_{f;\delta}(H,N_1,\cdots, N_K)$
can be written as the form (5.7) with $N^{2/9}\ll m_1\ll N^{4/9}.$

\section{\bf Proof of Theorem 3}

Suppose $k\geq 2$ and $\ell\geq 1$ are two fixed integers, and $f\in \{\tau_k, \Lambda_\ell, \mu_\ell,\omega_\ell\}.$ Recall that
$F(n)=\sum_{n=d^2m}g(d)f(m)$ with $g(d)\ll d^\varepsilon.$

Let $  y= x^{82/155}, $ then $x/y=x^{73/155}.$    Similar to (5.1) we have
\begin{eqnarray}
 \ \ \ \ \ \ E_F(x)
 && =-\sum_{x^{73/155}<n\leq x^{82/155}}F(n)\psi\left(\frac xn\right)+ \sum_{x^{73/155}<n\leq x^{82/155}}F(n)\psi\left(\frac{x}{n+1}\right)
 \\&&\ \ \ \ \ \ \ \ \ \ \ \ \ \
  + O(x^{73/155+\varepsilon} ).\nonumber
\end{eqnarray}
It suffices  to bound the sum
\begin{eqnarray*}
R_{F,\delta}(N;x):=\sum_{n\sim N}F(n) \psi\left(\frac{x}{n+\delta}\right)\ \ (\delta=0,1)
\end{eqnarray*}
for  $x^{73/155}\ll N\ll x^{82/155}.$ By the expression $F(n)=\sum_{n=d^2m}g(d)f(m) $ we see that
$R_{F,\delta}(N;x)$ can be written as a sum of $O(\log x)$ expressions of the form
\begin{eqnarray*}
R_{F,\delta}(D, M;x):=\sum_{\stackrel{d\sim D, m\sim M}{d^2m \asymp N}}f(m)g(d) \psi\left(\frac{x}{d^2m+\delta}\right)\ \ (\delta=0,1),
\end{eqnarray*}
where $D\geq 1, M\geq 1, D^2M\asymp N.$ Taking  $  \mathcal{H}= MDx^{-73/155}$   in    Lemma 2.7 we get that
\begin{equation}
R_{F,\delta}(D, M;x)\ll |R_{F,\delta}^{*}(D, M;x) |+O(x^{73/155+\varepsilon}),
\end{equation}
where
\begin{eqnarray*}
R_{F,\delta}^{*}(D, M;x): =\sum_{1\leq h\leq \mathcal{H}} \alpha^*(h)
 \sum_{\stackrel{d\sim D, m\sim M}{d^2m \asymp N}}f(m)g(d) e\left(\frac{hx}{d^2m+\delta}\right)
\end{eqnarray*}
with $\alpha^*(h)\ll 1/h.$

By a splitting argument we get
\begin{eqnarray}
R_{F,\delta}^{*}(D, M;x) \ll \frac 1H\left| \sum_{h\sim H} c_h
\sum_{\stackrel{d\sim D, m\sim M}{d^2m\asymp N}}f(m)g(d) e\left(\frac{hx}{d^2m+\delta}\right)   \right|\log x
\end{eqnarray}
for some $1\ll H\ll \mathcal{H},$ where $c_h=\alpha^{*}(h)h\ll 1.$     So Theorem 3 follows from (6.1)-(6.3)  and the estimate
\begin{eqnarray}
\ \ \ \ \ S_{F,\delta}(H, D, M):= \sum_{h\sim H} c_h \sum_{\stackrel{d\sim D, m\sim M}{d^2m\asymp N}}f(m)g(d) e\left(\frac{hx}{d^2m+\delta}\right)    \ll Hx^{73/155+\varepsilon}.
\end{eqnarray}

We consider three cases.

{\bf Case I}.  $D\gg N^{1/6}$.

 In this case by trivial estimate we have
$$
 S_{F,\delta}(H,D,M)\ll HD^{1+\varepsilon}M \ll \frac{HN^{1+\varepsilon}}{D}\ll HN^{5/6+\varepsilon}\ll Hx^{41/93+\varepsilon}\ll Hx^{73/155}.
$$

{\bf Case II.}  $x^{1/155}\ll D\ll N^{1/6}.$

   Similar to   Lemma 5.2, we can show that
$S_{F,\delta}(H, D,M)$ can be written as a sum of type I sums of the form
$$S_I(H,D,M_1,M_2):=\sum_{h\sim H}c_h\sum_{\stackrel{d\sim D, m_1\sim M_1,m_2\sim M_2}{d^2m_1m_2 \asymp N}}g(d)a_{m_1}
e\left(\frac{hx}{d^2m_1m_2+\delta}\right)$$
with $M_1\ll (N/D^2)^{1/3}$ and $a_{m_1}\ll m_1^\varepsilon,  $ and type II sums of the form
$$S_{II}(H,D,M_1,M_2):=\sum_{h\sim H}c_h\sum_{\stackrel{d\sim D, m_1\sim M_1,m_2\sim M_2}{d^2m_1m_2 \asymp N}}g(d)a_{m_1}b_{m_2}
e\left(\frac{hx}{d^2m_1m_2+\delta}\right)$$
with $(N/D^2)^{1/3}\ll M_1\ll (N/D^2)^{1/2}$ and $a_{m_1}\ll m_1^\varepsilon, b_{m_2}\ll m_2^\varepsilon.$
 Since the proof is almost the same as that of Lemma 5.1, we omit the details.

For $S_I(H,D,M_1,M_2),$ we  estimate the sum over $m_2$ by Lemma 2 with the exponent pair $(2/7,4/7)$ and estimate the sums over
other variables trivially. We get that
\begin{eqnarray}
  S_I(H,D,M_1,M_2)
&&\ll \left(\frac{D^3M_1^2  M_2^2}{x} + H (Hx)^{2/7}D^{3/7} M_1^{5/7}\right)N^\varepsilon \\
&&\ll  \left(\frac{N^2}{Dx} +  \frac{H^{9/7}  x^{2/7}   N^{5/21}}{D^{1/21}}\right)N^\varepsilon\nonumber \\
&&\ll \left(\frac{N^2}{x} +   H^{9/7}  x^{2/7}   N^{5/21} \right)N^\varepsilon\nonumber\\
&&\ll Hx^{1394/3255+\varepsilon}\ll Hx^{73/155}\nonumber
\end{eqnarray}
by noting that $M_1 \ll (N/D^2)^{1/3}, D\gg 1, N\ll x^{82/155}, H\ll MDx^{-73/155}.$

For  $S_{II}(H,D,M_1,M_2),$ we get by (4.19) of Lemma  4.5 that
\begin{eqnarray}
\ \ \ \ \ \ \ \   S_{II}(H,D,M_1,M_2)x^{-\varepsilon}&& \ll  \frac{H^{\frac 12}N^{\frac 56}}{D^{\frac 76}}
+\frac{Hx^{\frac{1}{4}}N^{\frac{3 }{8}}}{D^{\frac{1 }{2}}}  +
\frac{Hx^{\frac{1}{6}}N^{\frac{7}{12}}}{D^{\frac{2}{3}}}+\frac{H^2x}{DN^2}
\\
&&\ll
  Hx^{73/155 } \nonumber
\end{eqnarray}
by recalling that $x^{73/155}\ll N\ll x^{82/155}, x^{1/155}\ll D\ll D^{1/6}, H\ll DMx^{-73/155}.$

From (6.5) and (6.6) we see that
(6.4) holds for   $x^{1/155}\ll D\ll N^{1/6}.$

{\bf Case III.}   $D\ll x^{1/155}.$

We have
\begin{eqnarray}
S_{F,\delta}(H, D,M)x^{-\varepsilon}&&\ll \sum_{d\sim D}\left|\sum_{h\sim H}c_h\sum_{\stackrel{  m\sim M}{d^2m \asymp N}}f(m)
e\left(\frac{hx}{d^2m+\delta}\right)\right|\\
&&\ll \sum_{d\sim D}\left|\sum_{h\sim H}c_h\sum_{m\sim N_d}f(m)
e\left(\frac{hx_d}{m+\delta_d}\right)\right|\nonumber
\end{eqnarray}
 where $x_d=x/d^2, N_d=N/d^2\asymp M, \delta_d=\delta/d^2.$
 We shall show that the estimate
\begin{equation}
 \sum_{h\sim H}c_h\sum_{m\sim N_d}f(m)
e\left(\frac{hx_d}{m+\delta_d}\right)\ll H x_d^{8/17+\varepsilon}
 \end{equation}
 holds.

 If $N_d\ll x_d^{8/17},$ then (6.8) is trivial. Now suppose $N_d\gg x_d^{8/17}.$
 It is easy to see that
 $H\ll N_dx_d^{-8/17}$ and $N_d\ll x_d^{9/17}.$
So (6.8) follows from the estimate (5.5).

From (6.7) and (6.8) we get the estimate
\begin{eqnarray*}
S_{F,\delta}(H, D,M)x^{-\varepsilon}\ll H\sum_{d\sim D} x_d^{8/17}\ll x^{8/17}D^{1/17}\ll Hx^{73/155}
\end{eqnarray*}
for $D\ll x^{1/155}.$

This completes the proof of Theorem 3.

\end{document}